\documentclass[12pt]{article}
\usepackage{mathpazo}
\usepackage{amsfonts,amssymb,amsmath}
\usepackage{accents}
\usepackage{enumerate}
\usepackage{epigraph}
\usepackage{color,soul}     
\usepackage[normalem]{ulem}  
\usepackage{graphicx}
\graphicspath{ {./images/} } 

\newtheorem{proposition}{Proposition}
\newtheorem{axiom}{Axiom}
\newtheorem{lemma}{Lemma}
\newtheorem{corollary}{Corollary}

\newtheorem{example}{Example}
\newtheorem{definition}{Definition}

\newtheorem{question}{Question}

\newenvironment{proof}[1][Proof]{\noindent\textbf{#1.} }{\ \rule{0.5em}{0.5em}}

\usepackage{xcolor}
\pagecolor{white}

\begin{document}

\title{\textbf{Axioms for Arbitrary Object Theory}\thanks{Thanks to Sam Roberts, Philip Welch, Giorgio Venturi, and Hannes Leitgeb for invaluable comments.}}

\author{Luca Steinkrauss and Leon Horsten \\ University of Konstanz}
\date{\today}
\maketitle


\begin{abstract}
\noindent We formulate and discuss a general axiomatic theory of arbitrary objects. This theory is expressed in a simple first-order language without modal operators, and it is governed by classical logic.
\end{abstract}



\section{Introduction}

In contemporary metaphysics, there are two main theories of arbitrary objects. First, there is the theory that was developed by Fine in  his seminal monograph on arbitrary objects (\cite[Chapter 1]{Fine1985}). and extended and applied by him in subsequent publications. We will call this the F-theory of arbitrary objects. Secondly, Horsten has developed a rival theory of arbitrary objects (\cite{Horsten2019}). We will call this the H-theory. 

Both F-theory and H-theory are expressed as metaphysical theories  in ordinary English, with some heavy hints towards putative fundamental principles governing arbitrary objects.
One main difference between these theories is that in F-theory, \textit{dependence} is a primitive relation, whereas in H-theory, it is a defined concept. Another important difference is that whereas in F-theory, for every predicate $\varphi$, there is a unique independent arbitrary $\varphi$, in H-theory, this is not the case. A third difference is that in F-theory, the principle of generic attribution plays a crucial role, whereas this is not the case in H-theory.

Fine strongly suspects that, nonetheless, H-theory and F-theory are in some sense equivalent.\footnote{See \cite[p.~616ff]{Fine2022}.} However, precise claims about senses of equivalence cannot  be definitively established until both H-theory and F-theory are formally expressed in model-theoretic or axiomatic terms. A good reason for developing a fundamental formal account of F-theory and of H-theory is that such accounts can then be drawn upon when these theories are put to use in philosophical applications.

The present article constitutes a modest contribution towards the goal of formalising arbitrary object theory. Specifically, we will develop and discuss an axiomatisation of H-theory, which we call AOT  (``Arbitrary Object Theory''). This theory AOT is formulated in a first-order logical language without modal operators, and is governed by classical logic. 

The theory AOT intends to be a fundamental and a fully general (and somewhat flexible) theory of arbitrary objects. Ideally, it intends to be a suitable formal framework for all legitimate applications of arbitrary object theory. It must therefore postulate a rich variety of arbitrary objects that are equipped with appropriate formal properties. We will see (section \ref{conclusion}) that in this paper we do not quite achieve this aim. However, we do maintain that the theory that is proposed here should play a foundational role for current and future applications of arbitrary object theory.

According to the proposed theory, arbitrary objects are organised in correlated systems, where each such system of arbitrary objects is abstracted from a system of particular objects. Although this makes AOT an abstractionist account of sorts, it is, from an ontological perspective, at bottom not a reductionist view, for arbitrary objects are treated as \emph{sui generis} entities. However, we will also sketch how a corresponding reductionist theory of arbitrary objects can be constructed.

The structure of this article is as follows. In  section \ref{genericattribution}, we show how a formalisation of F-theory must take place in a non-classical logic, and take this as a mark against F-theory. In section \ref{axioms}, we formulate and discuss the language and axioms of AOT, and explore some elementary properties of arbitrary objects according to AOT. In section 6, we turn to the corresponding reductionist account, before closing with some remarks about possible questions and directions for further research.


\section{Generic attribution}\label{genericattribution}


Arbitrary objects are entities that can be in various \textit{states} (or \textit{situations}), where they take different \textit{values}. If we ignore higher-order arbitrariness, then the values of arbitrary objects are \textit{particular} objects. Fine proposes a \emph{criterion of identity} for arbitrary objects \cite[p.~18]{Fine1985}, which takes the form of an extensionality principle. The details of this principle need not detain us.\footnote{For a discussion, see \cite[section 7.3.1]{Horsten2019}.} In addition, Fine takes a central law governing reasoning about arbitrary objects to be the \textbf{Principle of Generic Attribution}, which is the cornerstone of his theory of arbitrary objects ( \cite[p.~17]{Fine1985}):

\begin{definition}[PGA]
Let $a_G$ be an arbitrary $G$. Then for all generic\footnote{The distinction between generic and non-generic (aka \emph{classical}) predicates is a difficult one. Fortunately we can ignore it here: we will also work with predicates that according to Fine clearly have a generic reading.} predicates $\varphi$:
$$  \varphi (a_G) \Leftrightarrow \forall i: G(i) \rightarrow \varphi (i)     $$
\end{definition}

\noindent Implicit in PGA, there is a \emph{comprehension principle} for arbitrary objects: if $G$ is any (non-empty) predicate, then there \emph{exists} at least one arbitrary $G$.

At first glance, PGA seems very intuitive: after all, we expect the arbitrary odd number, say, to be odd itself, much like we expect the arbitrary man to be mortal. However, accepting PGA forces one to abandon classical logic, as is shown by the following argument:

\begin{proposition}\label{propagainstPGA}

Suppose that a totally arbitrary object exists, and at least one object exists. Then 
$$ PGA \vdash_{CL}\forall x: x=d ,  $$
where $d$ is any individual constant for an object, and where $\vdash_{CL}$ is derivability in classical first-order logic (CL).\newline

\begin{proof}

\noindent Consider the vacuous condition $x=x$, and call it $G(x)$. This is a generic condition that specifies a totally arbitrary object.

\noindent Consider also the condition $x=d$ (with $d$ any individual constant), and call it $\varphi (x)$.

\noindent Then reason as follows in classical logic:

\noindent 1. $\varphi(a_G) \vee \neg \varphi(a_G)$ \hspace{2,2cm} Principle of Excluded Third

\noindent 2. $(\forall x : \varphi (x)) \vee (\forall x :  \neg \varphi(x))$ \hspace{0.5cm} from 1, using PGA and CL

\noindent 3. $d=d$ \hspace{4cm} Law of self-identity

\noindent 4. $\forall x : x=d $\hspace{3,3cm} 2,3, CL
\end{proof}

\end{proposition}

\noindent Observe that the assumptions about the metaphysics of arbitrary objects in this proof are very modest.   In particular, Fine's distinction between dependent and independent arbitrary objects plays no role in this simple argument. Concerning principles of reasoning about arbitrary objects, it is a pure and simple argument from PGA: Fine's assumption of Extensionality for arbitrary objects plays no role. The ``comprehension assumption'' implicit in PGA does, of course, play a role in the argument.

The conclusion of the argument of Proposition \ref{propagainstPGA} is clearly unacceptable. Therefore the argument has to be blocked somehow.

A somewhat more detailed but routine analysis\footnote{We leave this analysis to the reader.} of the reasoning in the proof of Proposition \ref{propagainstPGA} shows that in the transition from step 1 to step 2 we have implicitly used the following natural deduction rule of \textbf{Disjunctive Syllogism}:

\begin{center}

$\phi \vee \psi$

$\phi \rightarrow \tau$

$\psi \rightarrow \tau$

---------

$\tau $

\end{center}

Unlike in classical first-order logic, the rule of Disjunctive Syllogism is not valid in \emph{supervaluation logic}, which has the same theorems as the former, but fewer admissible inference rules. Fine indeed suggests that the way in which generic truth is treated in his theory can be subsumed under supervaluational truth.\footnote{See  \cite[p.~45]{Fine1985}.} We speculate  in this context that it is no coincidence that Fine is also the originator of the supervaluation theory of vagueness  \cite{Fine1975}: perhaps Fine's work on vagueness may have influenced his subsequent theory of arbitrary objects.
 
Withdrawing from classical logic is a heavy price to pay for salvaging a metaphysical theory: scientific theories are quite generally formulated in classical logic. Unlike classical logic, supervaluation logic is not recursively axiomatizable.\footnote{See \cite[p.~231]{KremerKremer2003}.} So even in principle we are only  able to reason in supervaluation logic to a limited extent. Therefore, if doing so does not incur too heavy a philosophical cost, retreating from classical logic to supervaluation logic should be avoided. We therefore now turn to our attempt at axiomatising arbitrary object theory in classical logic.


\section{Axioms}\label{axioms}

We now turn to the formulation and discussion of the axioms of AOT. Given any non-trivial collection of particular objects $P$, we want AOT to postulate a generous collection of (systems of) arbitrary objects that take values from $P$. In doing so, a \emph{right balance} needs to be struck between ontological generosity on the one hand and avoiding redundancy in the systems of arbitrary objects on the other hand.

The axiomatic theory AOT is formulated in a simple first-order language with identity. The language contains only two arbitrary object theory-specific non-logical symbols. Moreover, it does not contain modal operators;\footnote{In this respect, AOT differs from the theories that are proposed in \cite{HorstenSperanski2019} and in \cite{VenturiYago2024}.} Instead, it quantifies directly over possible situations or states (which are its counterparts of possible worlds).
In this way, the ontological commitments of AOT can directly be read off in Quinean fashion from the theory.

First-order classical logic is fully upheld. Given Proposition \ref{propagainstPGA}, this implies that Fine's principle PGA had better not be provable in AOT, and we will see that indeed it isn't.

We will further see that AOT contains a modest set theory, which plays a purely auxiliary role. Arbitrary objects, particular objects, and states of arbitrary objects are treated as \emph{Urelements} by this set theory.

AOT is a \emph{typed} theory of arbitrary objects in the sense that it maintains a strict separation of specific and arbitrary objects. Finally, as previously mentioned, it is a \emph{basic} theory of arbitrary objects in the sense that it does not deal with higher-order arbitrary objects. The task of formulating a formal theory of arbitrary objects that take other arbitrary objects (or perhaps even themselves) as values, is left for future work.


\subsection{Some basic elements of H-theory}\label{basic}

The H-theory of arbitrary objects is expressed informally in \cite{Horsten2019}. Let us start by briefly summarizing the basic elements of this view.

The man in the Clapham omnibus can be Laura Smith, living in 35 Holburn Avenue, Belford, but it can also be someone else. The man in the Clapham omnibus is an \emph{arbitrary} object; Laura is a \emph{particular} (or \emph{specific}) object.
Arbitrary objects are thus generally objects that can be in a state of 'being' (or, more precisely, taking as their value) one of many particular objects. Particular objects, in turn, are not the sort of objects that can be in such states. This indicates that basic conceptual components of arbitrary object theory are the notions of \emph{arbitrary object}, \emph{specific object}, \emph{taking a value}, and \emph{state}.

We are familiar from physics with the notion of a state. Classical mechanics, for instance, speaks of a state that a three particle system can be in.
Provisionally, a state can thus be regarded simply as a \emph{situation}, in the sense of \cite{BarwisePerry1983}: as an \emph{incomplete Leibnizian possible world}. For instance, there is the state, or situation, where Stephen is baking bread, and there is no fact of the matter whether in this state, Alice is working in the garden. Likewise, for any state of our three particle system, there is no fact of the matter whether, in this state, the man on the Clapham omnibus is Laura Smith.

Now, among the complete possible worlds, there is one unique \emph{actual} world. However, for states that arbitrary objects can be in, the notion of actuality does not make sense. It is meaningless to ask which state the arbitrary natural number is \emph{actually} in, or who the man on the Clapham omnibus \emph{actually} is: none of the possible situations that arbitrary objects can be in is metaphysically privileged in this way.\footnote{For this reason, \cite{Horsten2019} introduces a sui generis modality to describe these situations, called \emph{aftheiretic} possibility.}



\subsection{Language and definitions}

The language $\mathcal{L}_{AOT}$ in which AOT is formulated is a first-order language with identity, to which the following non-logical symbols are added:
$$ \in, Val, F .$$
$\in$ is a two-place predicate that expresses the \emph{elementhood} relation of the built-in set theory that plays an auxiliary role. The predicate $Val$ expresses a three-place relation. $Val(x,y,z)$ intuitively means that the value of arbitrary object $x$ in state $y$ is the particular object $z$. $F$ is a two-place predicate; we leave the discussion of this symbol until later.


In terms of the primitives $\in$ and $Val$, the notions of being a \textbf{particular object}, being an \textbf{arbitrary object}, being a \textbf{state}, and the notion of being a \textbf{set} can be defined as follows:

\vspace{0.2cm}

\noindent ``$z$ is a particular object'':
\begin{definition}\label{particularobject}
$P(z) \equiv \exists x \exists y Val(x,y,z)$
\end{definition}

\vspace{0.2cm}

\noindent ``$x$ is an arbitrary object'':
\begin{definition}\label{arbitraryobject}
$A(x) \equiv \exists y \exists z Val(x,y,z)$
\end{definition}

\vspace{0.2cm}

\noindent ``$y$ is a state'':
\begin{definition}\label{state}
$S(y) \equiv \exists x \exists z Val(x,y,z)$
\end{definition}

\noindent As we have seen in Section \ref{basic}, there is no need to be able to single one or more of the states as being somehow actual.

The particular objects, the arbitrary objects, and the sets together form the Ur-elements of the built-in set theory. Thus we define the predicate of being a set as follows:


\noindent ``$x$ is a set'':
\begin{definition}
$Set(x) \equiv \neg P(x) \wedge \neg A(x) \wedge \neg S(x)$
\end{definition}

\noindent So sets are not among the entities that arbitrary objects can take as values.

Since we will make sure that the sets, the particular objects, the arbitrary objects, and the states form mutually exclusive collections,\footnote{See Axiom \ref{disjoint} below.} we can and will use the notation of many-sorted quantification. In the interest of economy of expression, we will use restricted variables as follows:
\begin{itemize}
\item $a,b,c,\ldots$: arbitrary objects
\item $p,q,r,\ldots$: particular objects
\item $s,t,u,\ldots$: states
\item $m,n,o,\ldots$: sets
\end{itemize}

\noindent Thus, for example:
$$\forall a \phi(a) \equiv \forall x(A(x) \rightarrow \phi(x)),$$
and the obvious analogous conventions hold for the other restricted variables.


\subsection{Background axioms}

In what follows, we  sometimes we use some set theoretic or informal notation for concepts that can easily (but sometimes ponderously) be expressed in $\mathcal{L}_{AOT}$.

The logical principles of AOT are the usual axioms and rules of classical first-order logic with identity.
To this, we add axioms stating that the laws of a theory of sets with Urelements govern the predicates $Set$ and $\in$. There is no need here to be specific about the precise Urelement-set theory that we are adopting here, save to say that we  include the Axiom of Separation (for the whole language $\mathcal{L}_{AOT}$):

$$ \forall n \exists m \forall x [ x \in m \leftrightarrow (x \in n \wedge \phi (x))] \textrm{ for all } \phi \in  \mathcal{L}_{AOT} .  $$

\noindent In addition, we \emph{may} want to postulate that $P$ and $S$ are set-sized. It will turn out that if we do this, then systems of arbitrary objects will likewise be set-sized, and even the collection of all arbitrary objects will be set-sized. So this in effect amounts to adopting a \emph{set Urelement axiom}.\footnote{For a discussion of the set Urelement axiom, see \cite{McGee1997}.}


Now we discuss basic axioms that govern the proper AOT-predicate $Val$.  We start with very basic ones.

\begin{axiom}\label{oneparticularobject}
There is at least one particular object.
\end{axiom}
The reason for this is that arbitrary objects take particular objects as values. If there are no particular objects, then there are no arbitrary objects, either.

\begin{figure}[h]
\centering
\includegraphics[width=.5\textwidth]{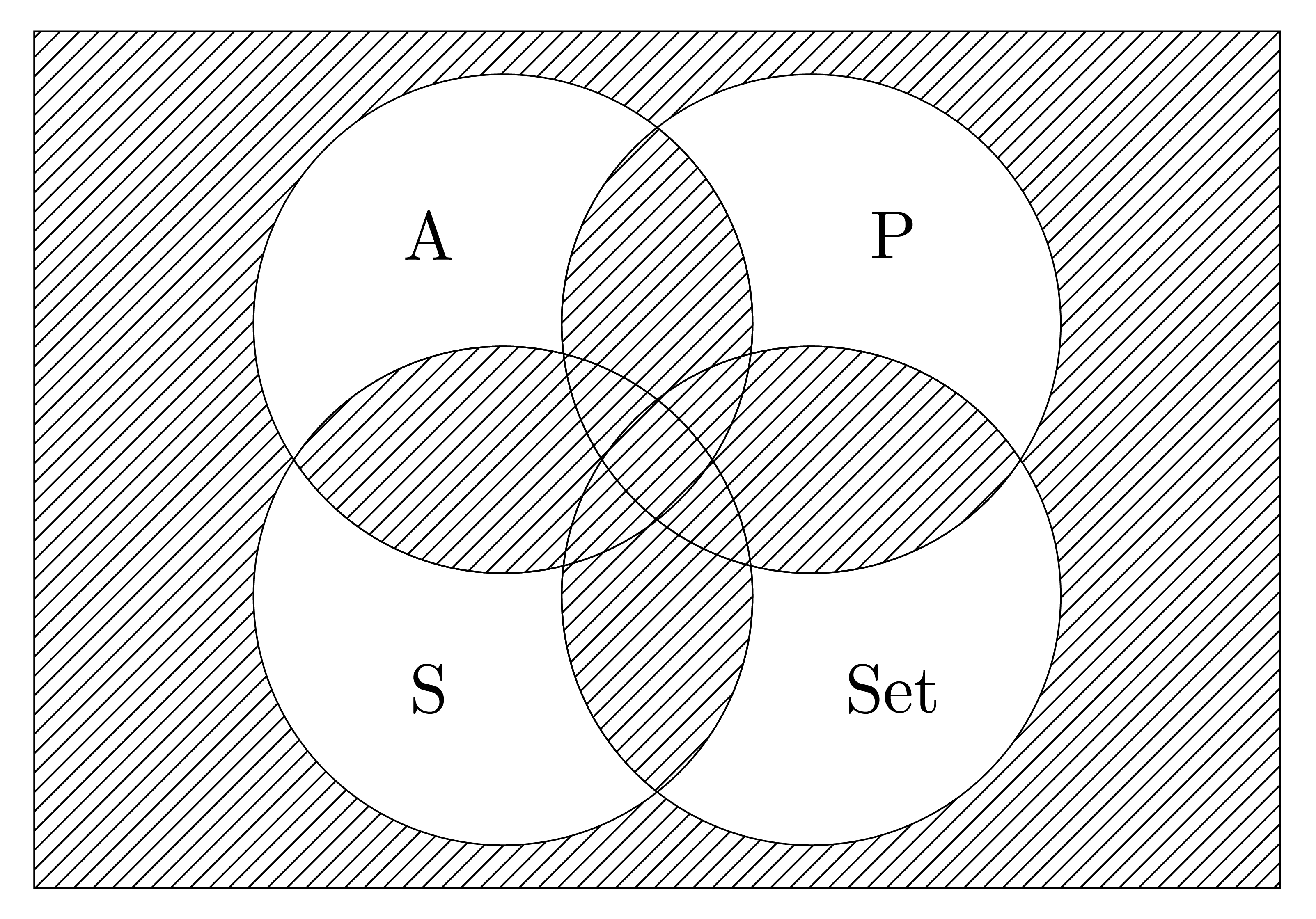}
\caption{The partition of the domain of discourse by $P,S,A,Set$. (The hatched area is empty.)}
\label{fig:disjoint}
\end{figure}

\begin{axiom}\label{disjoint}
The extensions of $P, A, S $ are pairwise disjoint.
\end{axiom}
Given our definition of $Set$, this means that we have \emph{partitioned} our domain of discourse into particular objects, arbitrary objects, states, and sets (see fig. \ref{fig:disjoint}). Thus, for instance, no entity can be both a state and an arbitrary object. From Axiom \ref{disjoint} it also follows immediately that GPA, as anticipated, does not hold in AOT.

\begin{axiom}
$Val$ is a partial function on $A\times S$.
\end{axiom}
This means that values of arbitrary objects in states are unique, while enabling the use functional notation when talking arbitrary objects (i.e., $Val(a,s)=p$ instead of $Val(a,s,p)$). Observe that we do not require arbitrary objects to take values in \emph{all} states. Indeed, it will turn out that (except in a trivial case) AOT does not allow arbitrary objects to take values in all states.



Next, we stipulate that every arbitrary object can only be in set-many states:
\begin{axiom}\label{setstates}
For every $a\in A$, there are only set many $s\in S$ such that $\exists z Val(a,s,z)$.
\end{axiom}
\noindent Since an arbitrary object is uniquely determined by which values it takes in which states (in a way that will be explained in the next section), this means that individual arbitrary objects can be seen as ``set-like'', i.e., \emph{relatively small} entities.

Axiom \ref{setstates} allows us to define \textbf{state spaces} (i.e. the set of states in which it takes values) and \textbf{value ranges} (i.e. the set of values it takes) as particular sets:

\noindent ``the state space of arbitrary object $a$'':
\begin{definition}
$ST(a)\equiv \{ y \mid \exists zVal(a,y,z)\}.$
\end{definition}

\noindent ``the value range of arbitrary object $a$'':
\begin{definition}
$VR(a)\equiv \{ z \mid \exists yVal(a,y,z) \}.$
\end{definition}


\subsection{Abstraction}

From Axiom \ref{oneparticularobject}, together with Definitions \ref{state} and \ref{arbitraryobject}, it already follows that there exist at least one arbitrary object and at least one state.  But even if we know that there are many particular objects, it is compatible with all that has been postulated so far that there is only one arbitrary object and only one particular object. So we want to state an axiom that guarantees that if there is a sizeable collection of particular objects, then there are also rich and varied collections of arbitrary objects taking values in rich collections of states. 

Such an axiom can be seen as a \textbf{comprehension axiom} for arbitrary objects. We have seen in section \ref{genericattribution} how in F-theory, comprehension was not fully articulated but mostly left implicit. 
In this section, we will see how formulating a suitable comprehension axiom for arbitrary objects in the context of classical logic is a somewhat delicate undertaking.

We start with defining what it means to be a \emph{particular object system}. We will then explain how from a particular object system, a \emph{system of arbitrary objects} can be abstracted.

\begin{definition}\label{pos}
For any ordinal $\lambda$, a \textbf{ particular object system} $o$ is any set of sequences of uniform length $\lambda$, where the elements of the sequences are particular objects.  
\end{definition}
\noindent Let $\mathcal{P}$ be the class of particular object systems.


A particular object system of $\lambda$-tuples of particular objects can be seen as an $\alpha \times \lambda$ matrix $M$---modulo permutations of the rows of $M$,--- for some ordinal $\alpha$. Informally, the columns $c_{\beta}$ of $M$  can be regarded as \emph{arbitrary objects}, and the rows $r_{\gamma}$ of $M$ can be regarded as the \emph{states} that these arbitrary objects can be in, where $c_{\beta}$ takes the $\beta$-th element of $r_{\gamma}$ as its value in ``state'' $r_{\gamma}$ (see fig. \ref{fig:abstraction} below).


Earlier we stated that set theory (with Urelements) merely serves us as a \emph{tool}, and that the sets, the arbitrary objects, and the states, are non-overlapping and therefore \emph{sui generis} ontological categories. Therefore it would be incorrect to say that arbitrary objects \emph{are} columns in particular object systems. Instead, particular object systems are in the final analysis only \emph{blueprints} from which systems of correlated arbitrary objects are \emph{abstracted}. In the expression of the abstraction relation, the function $F$ comes into play.

The next axiom describes how $F$ maps elements of particular object systems, indexed by the particular object systems to which they belong, onto the states.
\begin{axiom}\label{axiomF}
$F(x,y)$ is a two-place function such that:
\begin{enumerate}
\item (Domain) $F$ is defined on $x,y$ if and only if $x$ is a particular object system and $y\in x$;
\item the range of $F$ is $S$;
\end{enumerate}
\end{axiom}




Since we will often be interested in how an arbitrary object system is abstracted from a single given particular object system $o$, we will sometimes write $F_{o}(x)$ for the one-place function $F(o , x)$. As we will later see, the one-place function $F_{o}(x)$, which is defined only on $o$, is one-to-one.\footnote{It is not hard to see that if AOT were expressed in a second-order language, then there would be no need for a primitive symbol `$F$' in the language: the function $F$ would then just be definable. Since we have taken the conservative line of formalising arbitrary object theory in a first-order setting, this course of action is not open to us.}


Now we turn to the \textbf{comprehension principle} for arbitrary objects. In terms of the notion of particular object system and the function $F$, we can express how systems of arbitrary objects are \emph{abstracted} from systems of particular objects:

\begin{axiom}\label{abstractionaxiom}
For any particular object system $o$ consisting of tuples of length $\lambda$, there is a unique sequence $\mathcal{A}(o) = \langle a_1,\ldots,a_{\lambda}\rangle$ of arbitrary objects such that:
\begin{enumerate}
\item for every $x \in o$ and $a_{\alpha}$ (with $1 \leq\alpha \leq \lambda$),
$$ Val( a_{\alpha}, F_o(x) ) = x_{\alpha},   $$
where $x_{\alpha}$ is the $\alpha$-th element of $x$.
\item for any state which is not a value of $F_o$, the values for all $a \in \mathcal{A}(o)$ under $F_o$ are undefined. Here (and in what follows) $\mathcal{A}^*(o)$ is the \emph{set} of the elements in the ordered sequence $\mathcal{A}(o)$.
\item for all arbitrary objects $b \in A \setminus \mathcal{A}^*(o)$ and for all $x\in o$, $Val(b,F_o(x))$ is undefined. 
\end{enumerate}
\end{axiom}
\noindent Thus  $\mathcal{A}$ is an \emph{operator} on the class of particular object systems, and $\mathcal{A}(o) $ is the \textbf{arbitrary object system} ``abstracted'' by the operator $\mathcal{A}$ from  $o$.\footnote{Axiom \ref{abstractionaxiom} postulates the existence a variety of what are called \emph{arbitrary object spaces}: see \cite[Section 3.8]{Horsten2019}. We observe here that AOT allows ``incomplete arbitrary objects spaces'' (\cite[p.~55]{Horsten2019}), since it is allowed by Axiom \ref{abstractionaxiom} that $\lambda < \lvert P \rvert$. (Incomplete arbitrary object spaces are not emphasised much in \cite{Horsten2019}.)}

Since in the definition of particular object systems we allow the length $\lambda$ of sequences to be transfinite, AOT allows for infinite arbitrary object systems. For many applications of arbitrary object theory, only finite systems are of course needed. This just means restricting the value of $\lambda$ to values below $\omega$.

\begin{figure}[h]
\includegraphics[width=1\textwidth]{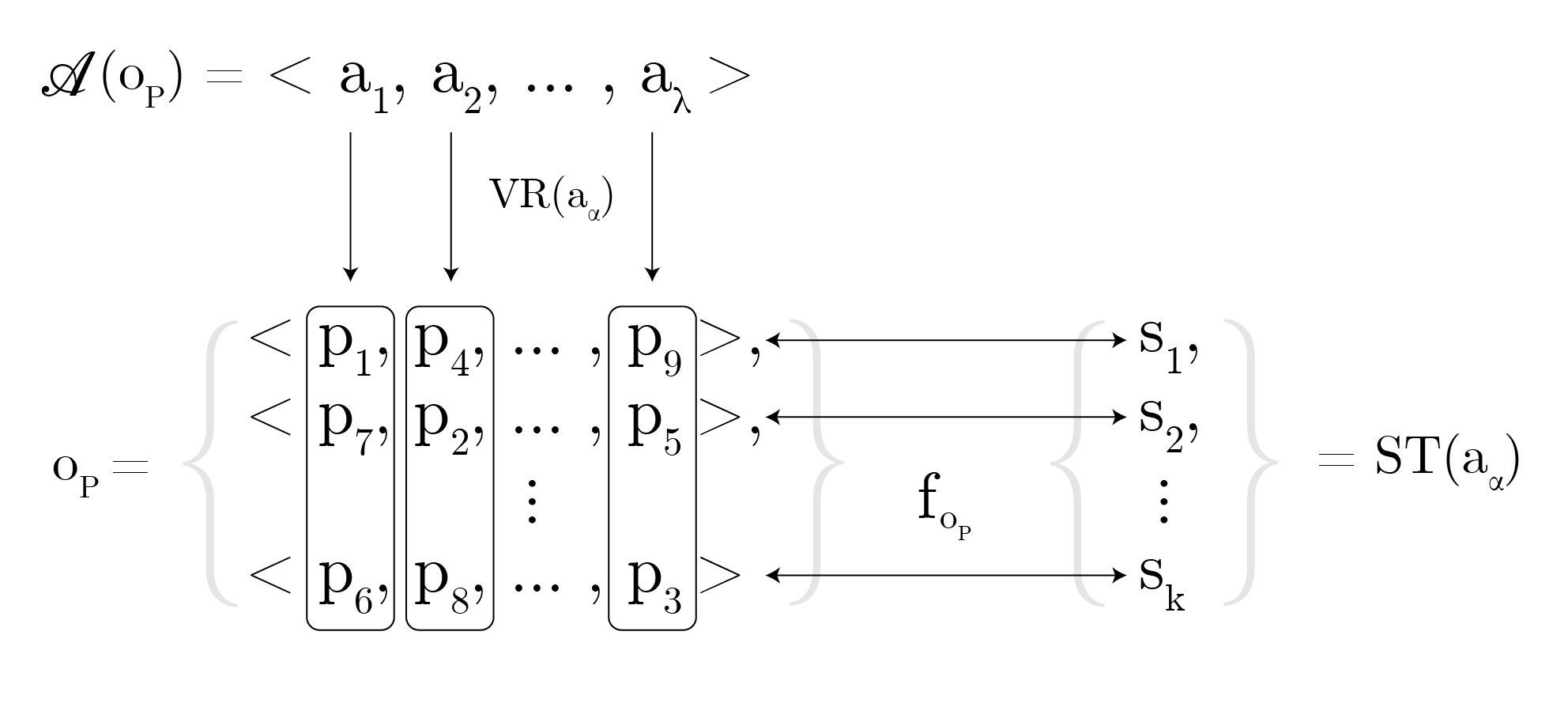}
\caption{Illustration of the structure postulated by Axiom \ref{abstractionaxiom}}
\label{fig:abstraction}
\end{figure}

We see that axiom \ref{abstractionaxiom} entails:
\begin{proposition}
For each particular object system $o$, and for each $a_{\alpha}, a_{\beta} \in \mathcal{A}(o): ST(a) = ST(b)$.
\end{proposition}

So we can unambiguously associate state spaces also with arbitrary object systems:
\begin{definition}
The state space of an arbitrary object system $\mathcal{A}(o)$ (denoted as $ST(\mathcal{A}(o))$ is $ST(a)$, for any $a\in \mathcal{A}(o)$.
\end{definition}
\noindent Observe that we may equivalently speak of the state space of the underlying set $\mathcal{A}^*(o)$ of the ordered sequence $\mathcal{A}(o)$. Indeed, we will in what follows often move quite freely between the ordered sequence $\mathcal{A}(o)$ and its underlying set $\mathcal{A}^*(o)$.

We next want to add the closure axiom which states that \emph{all} arbitrary objects belong to arbitrary object systems abstracted from particular object systems:
\begin{axiom}\label{closureaxiom}
Any arbitrary object belongs to $\mathcal{A}(o)$ for some particular object system $o$.
\end{axiom}
\noindent This can be seen as a kind of closure postulate for arbitrary objects.

Concluding the description of the axioms of AOT, we turn to the problem of \emph{criteria of identity} for arbitrary objects, which can be seen as a kind of \textbf{extensionality} question for arbitrary objects. The idea is to avoid needless duplication of abstracted entities. We distinguish between internal and external extensionality.

First,  an \textbf{internal} extensionality axiom prohibits unnecessary duplication of arbitrary objects \emph{within} arbitrary object systems. It does so by stipulating that arbitrary objects within one single arbitrary object system that have \emph{identical modal profiles}, are numerically identical:

\begin{axiom}\label{internalextensionality}
Let $o$ be any particular object system. Take any $a_{\alpha}, a_{\beta}\in \mathcal{A}^*(o)$ such that for each $s\in ST(\mathcal{A}(o)) $: $$ Val (a_{\alpha}, s)  = Val(a_{\beta}, s) .$$ Then $a_{\alpha} = a_{\beta}.$
\end{axiom}
\noindent In other words, any two identical columns in a particular object system determine one and the same abstracted arbitrary object.

Second, an \textbf{external} extensionality axiom prohibits needless duplication of arbitrary object systems. This axiom prohibits the existence of distinct arbitrary object systems that can be obtained from each other by a re-labelling of states or arbitrary objects:
\begin{axiom}\label{externalextensionality}
For any arbitrary object systems $A_1, A_2$, if there is a bijection f between $ST(A_1)$ and $ST(A_2)$ such that for all $a \in A_1$, there is an $a' \in A_2$ such that for all $s\in ST(A_1)$, $$ \forall z Val(a,s,z) \leftrightarrow Val(a',f(s),z),$$ and for all $a\in A_2$, there is an $a' \in A_1$ such that or all $s\in ST(A_2)$, $$ \forall z Val(a,s,z) \leftrightarrow Val(a',f^{-1}(s),z),$$
then $A_1 =A_2$.
\end{axiom}

\begin{figure}[h!]
\centering
\includegraphics[width=0.5\textwidth]{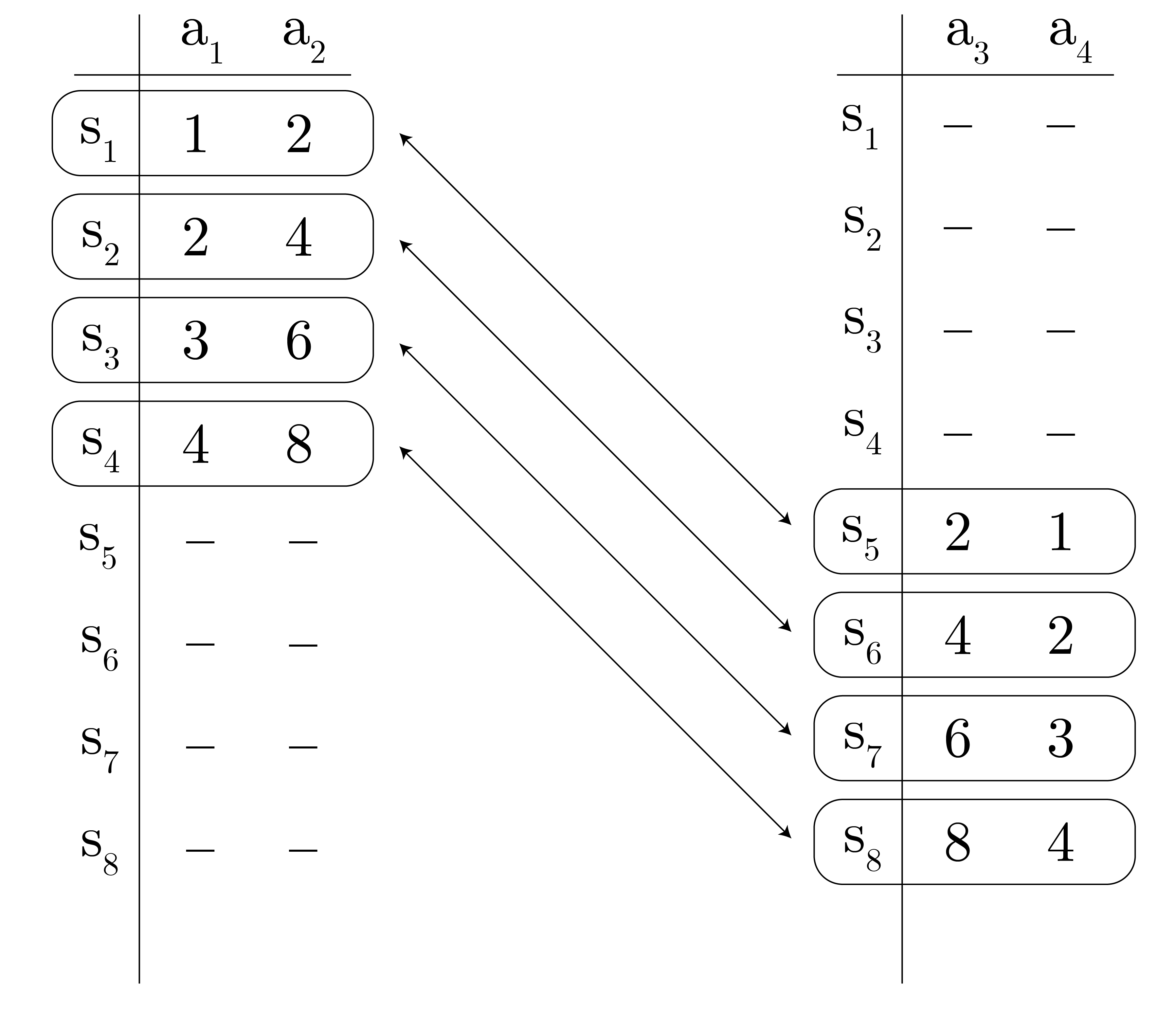}
\caption{Illustration of Axiom 10.}
\label{fig:identification}
\end{figure}

Axiom \ref{externalextensionality} prevents situations where two (or more) arbitrary object system have the same \emph{modal profile}  from occurring, by identifying them with each other (see fig. \ref{fig:identification}).
\\
\\
As a result, Axiom \ref{externalextensionality} forces identifications of states that would otherwise prevent the two-place function $F$ from being one-to-one:
\begin{proposition}
$F(x,y)$ expresses a 1-1 function.
\end{proposition}


\section{Simple examples}

The following two abstract but simple examples illustrate how arbitrary object systems are determined by particular object systems:
\begin{example}
Let $Q = \{ p_1, p_2, p_3  \}$, and consider the particular object system $o_Q = \{ \langle p_1, p_2, \rangle , \langle p_2, p_3 \rangle  \}$ with domain $Q$. Then $\lambda = 2$, and Axiom \ref{axiomF} postulates that there are states $s_1, s_2$ such that 
\begin{itemize}
\item $F_{o_Q}( \langle p_1, p_2, \rangle) = s_1$;
\item $F_{o_Q}(\langle p_2, p_3 \rangle)= s_2$.
\end{itemize}
Moreover, Axiom \ref{abstractionaxiom} postulates that there are arbitrary objects $a_1, a_2$ such that $\mathcal{A}(o_Q) = \langle a_1, a_2 \rangle$, where:
\begin{enumerate}
\item $Val(a_1, s_1)= p_1$;
\item $Val(a_1, s_2)= p_2$;
\item $Val(a_2, s_1)= p_2$;
\item $Val(a_2, s_2)= p_3$.
\end{enumerate}
\noindent This can be represented in a matrix (see fig. \ref{fig:example1} below).
\end{example}

\begin{figure}[h]
\centering
\includegraphics[width=0.25\textwidth]{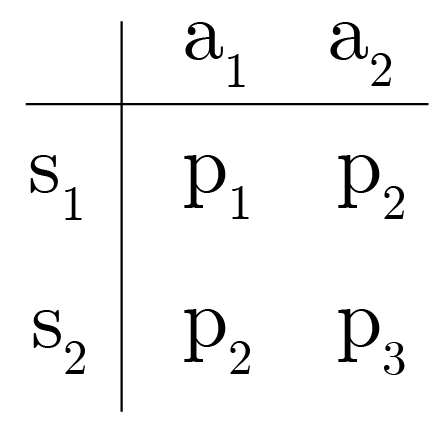}
\caption{Representation of Example 1.}
\label{fig:example1}
\end{figure}

\begin{example}
Suppose $P = \{0 ,1  \}$, and consider the particular object system $$o = \{\langle 0 ,1 \rangle , \langle 1, 0 \rangle  \}.$$ Then $\mathcal{A}(o)$ has a state space of two elements, and $\mathcal{A}(o)$ contains exactly two arbitrary objects $a_1$ and $a_2$. 

\noindent Consider the function $f$ such that $f(0)= 1$ and $f(1) = 0$. Then we see that for every state $s$, $Val(a_2, s, 0) \Leftrightarrow Val(a_1, s, 1)$, and $Val(a_2, s, 1) \Leftrightarrow Val(a_1, s, 0)$. This means that for every state $s$, the value of $a_2$ is determined by the value of $a_1$ in $s$ (and vice versa).
\end{example}

This elicits a general notion of \textbf{dependence} of arbitrary objects on other arbitrary objects:\footnote{Cfr \cite[section 9.4]{Horsten2019}.}

\begin{definition}
Let $a, b$ belong to some arbitrary object system $A$. Then we say that $b$ \textbf{depends on} $a$ iff there is a function $f$ on $P$ such that for every $s \in ST(\mathcal{A})$: 
$$Val(b,s, f(p)) \Leftrightarrow Val(a,s,p)$$ for any $p\in P$.
\end{definition}
\noindent In this way, we can capture something akin to a Finean notion of dependence between arbitrary objects\footnote{See \cite[p.~17]{Fine1985}.} in the present framework. Note, however, that, unlike in F-theory, in AOT, if $b$ depends on $a$, this does not entail any kind of metaphysical primacy of the latter over the former. Indeed, since dependence in AOT is nothing over and above a correlation of values in states, if the inverse of the dependence-defining relation is itself functional, then $a$ can be said to depend on $b$ as much as $b$ depends on $a$.

Let us now turn to a different (and somewhat more complicated) example:
\begin{example}
Suppose that our background theory of particular objects is a theory of the natural numbers, considered as particular objects. Then consider the particular object system $\{ \langle n \rangle \mid  n \in \mathbb{N}\}$, and the arbitrary object system $\mathcal{A}(\{ \langle n \rangle \mid n \in \mathbb{N}\})$ that it determines.
\end{example}
The arbitrary object system $\mathcal{A}(\{ \langle n \rangle \mid n \in \mathbb{N}\})$ contains exactly one arbitrary object $a_0$, namely the arbitrary object that can take each particular natural number as one of its values (and can take nothing else as its value). This arbitrary number $a_0$ can be seen as the ``Ur-arbitrary natural number'', since it depends on no other arbitrary natural numbers and in turn no other arbitrary natural numbers depend on it: it is the AOT-pendant of what Fine calls the ``\emph{initial} arbitrary natural number'' (\cite[p.~608]{Fine2022}).


\section{Basic properties}\label{properties}

Next, we turn to the discussion of some further elementary properties of arbitrary objects and arbitrary object systems according to AOT.  They are intended to give us a feeling of the nature of such entities according to AOT.\footnote{Where an argument for a proposition is bordering on triviality, we omit it.}

We first address questions concerning the existence and size of models of AOT.

\begin{proposition}
If there is exactly one particular object, then there exists one and only one arbitrary object.
\end{proposition}

Indeed, each ``degenerate'' one-element particular object system $$o_p= \{ \langle p \rangle \}$$ (for $p\in P$) is canonically embedded into $S$: each such $o_p$ determines exactly one arbitrary object $a$, which can only be in the state where it takes $p$ as its value.\footnote{Similarly, there are more complicated arbitrary object systems that contain ``arbitrary'' objects that in each state where they take a value, take the same particular value.}
If we do not want such degenerate arbitrary objects in our ontology, then it is of course easy to eliminate them by slightly modifying AOT; we will not do so here.

This inspires a simple  \emph{consistency proof}:
\begin{proposition}
AOT is consistent.

\begin{proof}
Consider the model in which $P = \{ p \}  $, $S =  \{ s \}$, $A =  \{ a \}$. Then $F( \{  \langle p \rangle \} , \langle  p \rangle) = s$, $Val(a,s) = p$, and it can be routinely verified that the axioms of AOT are satisfied.
\end{proof}
\end{proposition}

\begin{proposition}
If there are at least two particular objects, then there are infinitely many arbitrary object systems, and arbitrary object systems with unboundedly many arbitrary objects.

\begin{proof} 
Let 0 and 1, respectively, be the only particular objects in the domain of discourse. Consider the class of \emph{1-0 -diagonal} arbitrary object systems, which are those where the $\alpha$-th arbitrary object takes 1 as a value in the $\alpha$-th state, and 0 otherwise. Suppose for a reductio, that there is a largest such arbitrary object system, which contains $n$ arbitrary objects, and as many states. From definition \ref{pos} we know that there is a particular object system with $n$+1 sequences of length $n$+1, where the $\alpha$-th element in the $\alpha$-th sequence is 1, and every other element is 0. Due to Axiom \ref{abstractionaxiom}, there will be an arbitrary object system abstracted from the described particular object system, which is both 1-0-diagonal and has $n$+1 arbitrary objects, against our initial assumption.


\end{proof}
\end{proposition}

\begin{figure}[h]
\centering
\includegraphics[width=0.5\textwidth]{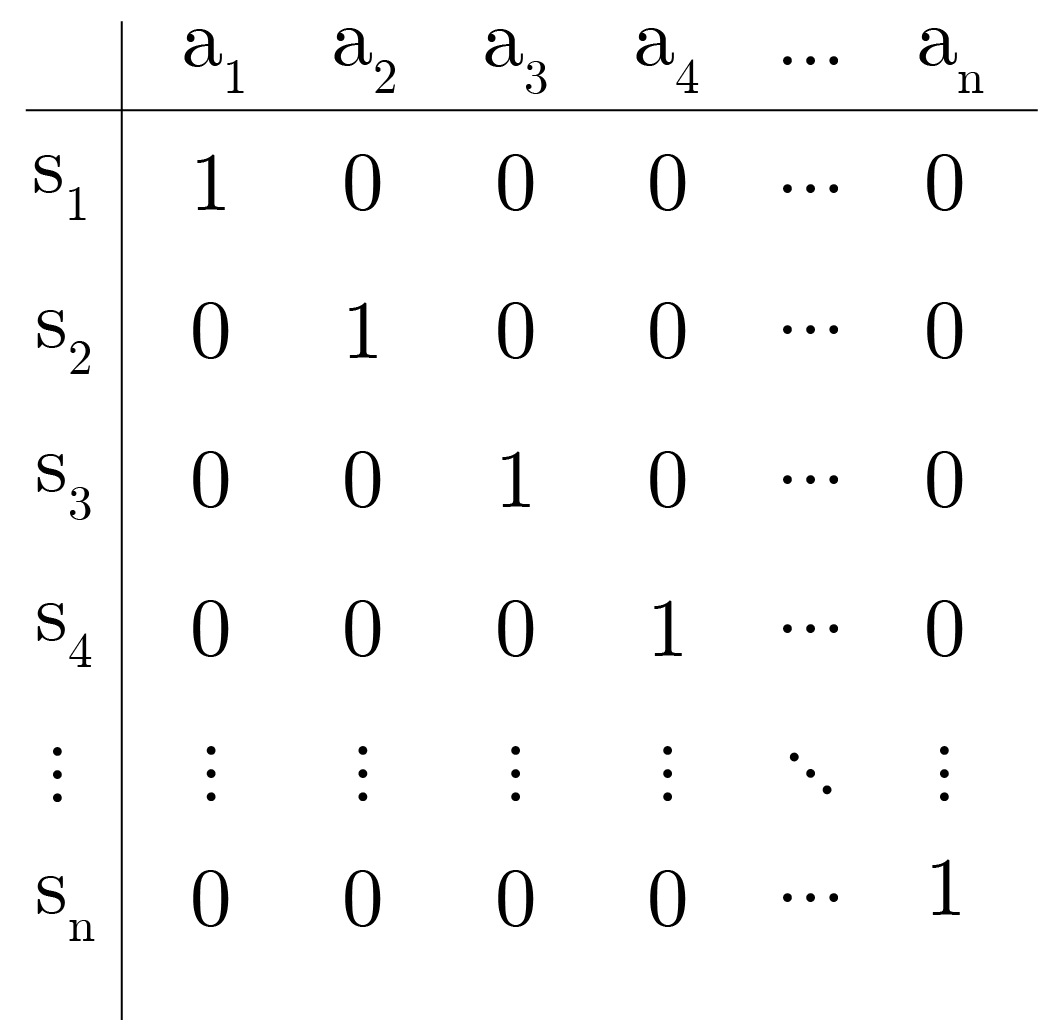}
\caption{A 1-0-diagonal arbitrary object system of size $n$.}
\label{fig:diagonal}
\end{figure}

\begin{proposition}
For all particular object systems $o$, and for each $a \in \mathcal{A}(o)$: $\lvert  o \rvert = \lvert  ST (a)  \rvert$.
\end{proposition}

\begin{definition}
For all particular object systems $o$, we define $l(o)$ to denote the length of the tuples of $o$.
\end{definition}

\begin{proposition}
For all particular object systems $o$: $\lvert \mathcal{A}(o) \rvert \leq \lvert l (o) \rvert$.
\end{proposition}
\noindent But observe that axiom \ref{internalextensionality} implies that there are particular object systems $o_P$ such that $\lvert \mathcal{A}(o_P) \rvert  <  \lvert l (o_P) \rvert$.

\begin{proposition}
Suppose that $\lvert P \rvert = n$. Then for any finite number $m$, there is a maximal finite number of states that any arbitrary object system $\mathcal{A}$ with $\lvert \mathcal{A} \rvert = m$ that $ST(\mathcal{A})$ can contain.

\begin{proof}
This follows from the definition of particular object systems as sets, rather than sequences) of sequences of particular objects, together with Axioms \ref{abstractionaxiom} and \ref{closureaxiom}. Indeed, it is easy to see that the maximal number of states will be $n^m$.
\end{proof}
\end{proposition}





Consider any two models $\mathcal{M}_1, \mathcal{M}_2$ of AOT that contain the same collection $P$ of particular objects. Given the fact that there is no categoricity theorem for the background (first-order) set theory of AOT, we can \emph{a fortiori} not in general expect that $\mathcal{M}_1$ and $\mathcal{M}_2$ are isomorphic. Nonetheless, given the fact that axiom \ref{externalextensionality} prevents $\mathcal{M}_1$ or $\mathcal{M}_2$ to have one or more isomorphic copies of arbitrary object systems, we can prove that the respective classes $\mathcal{C}_1, \mathcal{C}_2$ of \emph{finite} arbitrary object systems of $\mathcal{M}_1$ and $\mathcal{M}_2$---where $A$ being finite means not only containing only a finite number of arbitrary objects, but also the state space $ST(A)$ being finite---are isomorphic:
\begin{proposition}[limited categoricity]
For any two models models $\mathcal{M}_1, \mathcal{M}_2$ of AOT that contain the same collection $P$ of particular objects, their respective classes $\mathcal{C}_1, \mathcal{C}_2$ of finite arbitrary object systems are isomorphic.

\begin{proof}
Observe that due to axioms \ref{abstractionaxiom} and \ref{closureaxiom}, for any arbitrary object system belonging to some model of AOT, that model contains some particular object system it was abstracted from. Conversely, for every collection of particular object systems with isomorphic collapses belonging to some model of AOT, that model contains exactly one corresponding arbitrary object system, thanks to axioms \ref{internalextensionality} and \ref{externalextensionality}. Further, any two models of AOT $\mathcal{M}_1$ and $\mathcal{M}_2$ which contain the same particular objects also contain the same collection of finite particular object systems. Therefore, for any finite arbitrary object system $\mathcal{A}_1$ in $\mathcal{M}_1$, there must be exactly one arbitrary object system in $\mathcal{M}_2$, call it $\mathcal{A}_2$, which is abstracted from a particular object system whose collapse is isomorphic to the collapse of the particular object system that $\mathcal{A}_1$ is abstracted from, and vice versa. Since arbitrary object systems abstracted from particular object systems with isomorphic collapses are themselves isomorphic, any such arbitrary object systems $\mathcal{A}_1$ and $\mathcal{A}_2$ are isomorphic.
\end{proof}
\end{proposition}

At any rate, we conjecture that the number of arbitrary object systems containing $n$ arbitrary objects is described by a fast growing function (of $n$):
\begin{question}
Let $C_P(n)$ be the collection of all arbitrary object systems with at most $n$ arbitrary objects per system, where each arbitrary objects takes values from $P$. What is the cardinality of $C_P(n)$ as a function of $n$?\end{question}

Next, we discuss an important consequence from Axiom \ref{abstractionaxiom}:

\begin{lemma}\label{uniquestates}
For all particular object systems $o_1, o_2$:
$$ ST(\mathcal{A}(o_1)) \cap ST(\mathcal{A}(o_2)) = \emptyset \vee \mathcal{A}^*(o_1) = \mathcal{A}^*(o_2).$$

\begin{proof}
Suppose that the first disjunct is false, i.e., there is some $s \in ST(\mathcal{A}(o_1)) \cap ST(\mathcal{A}(o_2))$. Then there is an $x_1 \in o_1$ such that $F_{o_1}(x_1) =s$. Fix such an $x_1$. Similarly, there must be an $x_2 \in o_2$ such that $F_{o_2}(x_2) =s$. Fix such an $x_2$.

Now take any $a \in \mathcal{A}^*(o_1)$. Suppose, for a reductio, that $a \not \in \mathcal{A}^*(o_2)$, i.e., $a \in A \setminus \mathcal{A}^*(o_2)$. Then, by clause 3. of Axiom \ref{abstractionaxiom}, for all $y \in o_2$, $Val(a, F_{o_2}(y))$ is undefined. But $Val(a, F_{o_2}(x_2))$ is defined, and $x_2 \in o_2$. So we have reached a contradiction. Therefore $a \in \mathcal{A}^*(o_2)$.

By a symmetric argument, it can be seen that for every $a \in \mathcal{A}^*(o_2)$, $a$ must also be in $\mathcal{A}^*(o_1)$.
\end{proof}
\end{lemma}

Lemma \ref{uniquestates} says that any two arbitrary object systems are completely isolated from each other, in the sense that no arbitrary object belonging to one of them takes a value in any state in which an arbitrary object belonging to the other system takes a value.

If we think, for instance, of random variables as particular kinds of arbitrary objects,\footnote{See \cite[chapter 10]{Horsten2019}.}, then this makes sense. When two systems $S_1$ and $S_2$ of random variables are defined (on separate occasions, say), it is assumed that a variable in $S_1$ and does not depend on or correlate with any variable in $S_2$. Also in many other applications of arbitrary object theory, this assumption makes sense. However, it may be questioned if this assumption should hold for systems of arbitrary objects \emph{in general}. Suppose that we want to treat positions in structures as arbitrary objects belonging to arbitrary object systems.\footnote{See for example \cite[section VI]{Fine1998}.} According to some views, the natural number 1 is the same entity as the real number 1. If that is so, then the number 1, as an arbitrary object, would belong to two systems of arbitrary objects, the natural number structure and the real number structure. To accommodate such a view, one would have to omit clause 3. from Axiom \ref{abstractionaxiom}.

In combination with Axiom \ref{internalextensionality}, Lemma \ref{uniquestates} entails an identity criterion for arbitrary objects in general:

\begin{lemma} \label{identcrit}

$ \forall a,b [a=b \leftrightarrow \forall s \forall p (Val(a,s,p) \leftrightarrow Val(b,s,p))]$

\begin{proof}

\noindent The left-to-right direction of the equivalence trivially follows from Leibniz's law of the substitutivity salva veritate of identicals. So we focus on the right-to-left direction. 

\noindent Consider any $a,b$, and assume $\forall s \forall p (Val(a,s,p) \leftrightarrow Val(b,s,p)]$.

\noindent By Axiom \ref{closureaxiom}, there must be $o_1, o_2$ such that $a \in \mathcal{A}^*(o_{1})$ and $b \in \mathcal{A}^*(o_{2})$. Then by the law of excluded third:

\noindent (1) $ \mathcal{A}^*(o_{1}) = \mathcal{A}^*(o_{2})$. Then by Axiom \ref{internalextensionality} it follows that $a=b$.

\noindent (2) $ \mathcal{A}^*(o_{1}) \not = \mathcal{A}^*(o_{2})$. Then by Lemma \ref{uniquestates} it follows that $ST(\mathcal{A}(o_1)) \cap ST(\mathcal{A}(o_2)) = \emptyset $, whereby $Val(a,s,p) \rightarrow \neg Val(b,s,p)$. But this means that $\neg \forall s \forall p (Val(a,s,p) \leftrightarrow Val(b,s,p) )$, contrary to our assumption. So situation (2) cannot occur.

\noindent This means that only (1) is relevant, and the result follows.
\end{proof}
\end{lemma}

\begin{corollary}\label{cor}

For any particular object systems $o_1, o_2$:
$$  \mathcal{A}^*(o_{1}) = \mathcal{A}^*(o_{2}) \Rightarrow ST(\mathcal{A}^*(o_{1})) = ST (\mathcal{A}^*(o_{2})) .$$

\begin{proof}
Follows directly from Lemma \ref{identcrit}.
\end{proof}
\end{corollary}

Axiom \ref{internalextensionality} suggests the definition of the concept of the \emph{collapse} of a particular object system:
\begin{definition}
For any particular object systems $o$, let $col(o)$ (the collapse of $o$) be the result of removing from $o$ duplicates of columns, i.e., later columns for which $o$ contains an earlier column with identical entries.
\end{definition}

\begin{lemma}

For all particular object systems $o_1, o_2$:
$$   \mathcal{A}^*(o_{1})  = \mathcal{A}^*(o_{2}) \leftrightarrow col(o_1) = col(o_2).  $$

\begin{proof}
\noindent The right-to-left direction follows from Axiom \ref{internalextensionality}.

\noindent For the left-to-right direction, suppose, for a reductio, that $ \mathcal{A}^*(o_{1}) = \mathcal{A}^*(o_{2})$ but  $col(o_1) \not = col(o_2)$. 

\noindent From $ \mathcal{A}^*(o_{1}) = \mathcal{A}^*(o_{2})$, it follows by corollary \ref{cor} that $ST(\mathcal{A}^*(o_{1})) = ST (\mathcal{A}^*(o_{2}))$.

\noindent If $col(o_1) \not = col(o_2)$, then there is a column $c$ in $o_2$ that is not in $o_1$, or vice versa. (The latter is a symmetric case that leads to the same conclusion: we will not further discuss it here.) Then $c$ determines an arbitrary object $a_c$ such that for each $b \in \mathcal{A}^*(o_{2})$, it takes a different value in some state. So by Lemma \ref{identcrit}, we have $c \not \in \mathcal{A}^*(o_{2})$. So we have reached a contradiction with our assumption that $ \mathcal{A}^*(o_{1}) = \mathcal{A}^*(o_{2})$.
\end{proof}

\end{lemma}


\section{A definitional approach}\label{mathperspective}

Until now, we were motivated by metaphysical purity, and we eschewed any form of ontological reduction. But if we focus on the mathematical structure of AOT, then we see that its mathematical core is actually very \emph{simple}. In other words, the price of metaphysical purity is the existence of mathematically redundant structure.
This becomes clear  when we abandon the view that the states, arbitrary objects, and sets should be treated as mutually non-overlapping categories. In such a situation, the general theory of arbitrary objects can be formally presented and investigated in an Ur-element set theory, where the particular objects exhaust the Ur-elements.

Consider the first-order language $\mathcal{L}_{P, \in}$: it has the predicate $P$ and the elementhood relation $\in$ as its sole non-logical symbols. Moreover, let $U$ be some standard Ur-element set theory in $\mathcal{L}_{P, \in}$. The finer details of $U$ again do not matter for our discussion, but we will assume that the particular objects (i.e., the Ur-elements) form a non-empty \emph{set}.

We first define what it means to be a particular object system. This definition is as before, except that, for simplicity's sake, we make the simplifying assumption that the length of the tuples of any particular object system $o$ (denoted as $l(o)$) is \emph{finite}.
\begin{definition}
A particular object system $o$ is a set of tuples $\overline{p}$ of elements of $P$, where the tuples $\overline{p}$ are of uniform finite length $l$, and are such that for no two $m,n \leq l$, $\overline{p}(m) = \overline{p}(n)$ for each $\overline{p} \in o$.
\end{definition}

\noindent The finiteness of the length of the tuples in particular object systems, together with the set-sizedness of $P$, allows us to carry out the definition of the particular object systems in $\mathcal{L}_{P, \epsilon}$, using the separation axiom in $U$, as a particular impure \emph{set}. The ``no duplicate columns''-restriction on particular objects is imposed to ensure that one particular object system is not allowed to ``contain'' two duplicate arbitrary objects.

In terms of the notion of particular object system, we next \emph{define} the state predicate in $\mathcal{L}_{P, \in}$:
\begin{definition}
A state is an ordered pair $\langle o, \overline{p} \rangle$, with $o$ a particular object system and $\overline{p} \in o$.
\end{definition}
\noindent This means that, unlike in the approach that was discussed in the previous sections, states are now conceived of as impure sets. 

Next, we define in $\mathcal{L}_{P, \in}$ what it means to be an arbitrary object:
\begin{definition}
An arbitrary object is an ordered pair $\langle o, m  \rangle$, with $o$ a particular object system, and $m \leq l(o_P)$.
\end{definition}
\noindent For the same reasons as before, the collection of arbitrary objects is defined as a set in $\mathcal{L}_{P, \in}$, using the separation axiom.

To conclude, we can define the valuation function in in $\mathcal{L}_{P, \epsilon}$:
\begin{definition}
$Val(x,y,z)$ holds if $x$ is an arbitrary object $\langle o, m \rangle$, $y$ is a state $\langle o, \overline{p}\rangle$, $z$ is a particular object, and $\overline{p}(m) = z$. 
\end{definition}
So arbitrary objects can be seen as \emph{projection operators}---such as operators on Hilbert spaces, for instance---on states of particular object systems. Again, the valuation relation is then just a definable set in $\mathcal{L}_{P, \in}$.

Then it is easy to see that properties similar to those discussed in the previous section can be proved in $U$ about particular objects, states, arbitrary objects, and the valuation relation. In particular, the axioms of AOT will be verified by these models, which can be seen as the \emph{intended} models from a mathematical point of view.


\section{In conclusion}\label{conclusion}

In this article, we have proposed a basic general theory of arbitrary objects (AOT). First, we developed  this theory from a purely formal metaphysical perspective (sections \ref{axioms} and \ref{properties}); afterwards, we zoomed in on the purely mathematical core of the theory (section \ref{mathperspective}).

The fundamental features of the theory that we propose, such as the claim that arbitrary objects are organised in \emph{systems}, are essential to our proposal. Nonetheless, AOT is also intended to be a \emph{flexible} theory: it allows for variation. For instance, it was postulated that arbitrary object systems are isolated from each other: they share no states. However, the axiom that postulates this (Axiom \ref{uniquestates}) may be omitted from AOT without rejecting the core of the theory, thus allowing systems of arbitrary objects to overlap.

The formal theory of arbitrary objects AOT is intended to serve as a cornerstone for all legitimate applications of arbitrary object theory. It should provide fundamental elements for the formal frameworks for viewing structures as arbitrary entities, for viewing boolean-valued models as arbitrary objects, for viewing random variables as arbitrary objects\ldots

However, for several reasons, in this article we have not quite reached this aim. Firstly, in AOT the sets are kept disjoint from the particular objects. As a consequence, AOT does not contain a theory of arbitrary \emph{sets}. Similarly, the states are kept disjoint from the particular objects. Thus AOT also fails to accommodate arbitrary states. Secondly, AOT only describes first-order arbitrariness. Higher-order arbitrary objects, can take arbitrary objects (of lower order) as their value, are outside the scope of the theory.\footnote{Higher-order arbitrariness is involved in interpreting Boolean-valued sets as arbitrary objects: see \cite{Horsten2024}.} Lastly, AOT is a \emph{typed} theory. This entails that no arbitrary object can take \emph{itself} as a value. Thus also the development of a \emph{type-free} theory of arbitrary objects is left for future research. 
We express the hope that AOT can be adapted and generalised in natural ways so as to address these challenges successfully.


\end{document}